\documentclass[a4paper,12pt]{article}

\usepackage{fancyhdr}
\usepackage{epsfig}
\usepackage{cite}
\usepackage{amsmath}
\usepackage{theorem}
\usepackage{graphicx}
\usepackage{amssymb}
\usepackage{latexsym}
\usepackage{epic}
\usepackage{amsfonts,amssymb,amsopn}
\usepackage[all,2cell,dvips]{xy} \UseAllTwocells \SilentMatrices

\newcommand{\C}{\mathbb{C}}
\newcommand{\Mt}{{\cal M}_{2}}
\newcommand{\Mn}{{\cal M}_{n}}
\newcommand{\PP}{\mathbb{P}}
\newcommand{\Spn}{\mathrm{Sp}_{n}\C}

\newcommand{\Homom}{\mathrm{Hom}}

\newcommand{\Ker}{\mathrm{Ker}}
\newcommand{\degree}{\mathrm{deg}}
\newcommand{\rank}{\mathrm{rk}}
\newcommand{\Sym}{\mathrm{Sym}}
\newcommand{\strs}{{\cal O}_{X}}
\newcommand{\Image}{\mathrm{Im}}
\newcommand{\urd}{{\cal U}_{X}(r,d)}
\newcommand{\utmo}{{\cal U}_{X}(2,-1)}
\newcommand{\uto}{{\cal U}_{X}(2,1)}
\newcommand{\surL}{{\cal SU}_{X}(r,L)}
\newcommand{\cans}{{\cal K}_{X}}

\newcommand{\sRat}{\underline{\mathrm{Rat}}}
\newcommand{\Rat}{\mathrm{Rat}}
\newcommand{\sPrin}{\underline{\mathrm{Prin}}}
\newcommand{\Prin}{\mathrm{Prin}}
\newcommand{\Supp}{\mathrm{Supp}}

\newcommand{\qed}{\ifhmode\unskip\nobreak\fi\quad\ensuremath\square}

{\theorembodyfont{\slshape}\newtheorem{prop}{{\textbf Proposition}}}
{\theorembodyfont{\slshape}\newtheorem{thm}[prop]{{\textbf Theorem}}}
{\theorembodyfont{\slshape}\newtheorem{cor}[prop]{{\textbf Corollary}}}
{\theorembodyfont{\slshape}\newtheorem{lemma}[prop]{{\textbf Lemma}}}
{\theorembodyfont{\slshape}\newtheorem{crit}[prop]{{\textbf Criterion}}}

\title{Moduli of rank 4 symplectic vector bundles over a curve of genus 2}

\author{George H.\ Hitching}

\begin{document}

\maketitle

\begin{abstract}
Let $X$ be a complex projective curve which is smooth and irreducible of genus $2$. The moduli space $\Mt$ of semistable symplectic vector bundles of rank $4$ over $X$ is a variety of dimension $10$. After assembling some results on vector bundles of rank $2$ and odd degree over $X$, we construct a generically finite cover of $\Mt$ by a family of $5$-dimensional projective spaces, and outline some applications.
\end{abstract}

\tableofcontents

\section{Introduction}

Let $X$ be a complex projective smooth irreducible curve of genus $g$ with structure sheaf $\strs$ and canonical bundle $K_X$. In this section, we review the notion of a symplectic vector bundle over $X$ and introduce the moduli space we will be studying.\\

\subsection*{Symplectic vector bundles}

\textbf{Definition:} A \textsl{symplectic vector bundle} over $X$ is a pair $(W, \theta)$ where $W$ is a vector bundle over $X$ and $\theta$ is a bilinear nondegenerate antisymmetric form on $W \times W$ with values in a line bundle $L$.\\
\\
For us, $L$ will always be the trivial bundle $O_X$. Such a bundle $W$ is necessarily of even rank. If there is no ambiguity, we write $W$ for the pair $(W, \theta)$. For example, if $W$ is simple then $\theta$ is unique up to nonzero scalar multiple.
\par
A subbundle $E \subset W$ is \textsl{isotropic} if $\theta(E,E) = 0$. For any subbundle $E \subseteq W$, we have the short exact vector bundle sequence
\[ 0 \to E^{\perp} \to W \to E^{*} \to 0 \]
where the surjection is the map $w \mapsto \theta( w , \cdot )|_E$, and
\[ E^{\perp} = \left\{ w \in W : \theta( w , E ) = 0 \right\} \]
is the \textsl{orthogonal complement} of $E$. Clearly $E$ is isotropic if and only if $E \subseteq E^{\perp}$; this shows that the rank of an isotropic subbundle is at most $\frac{1}{2} \rank(W)$. An isotropic subbundle of maximal rank is called a \textsl{Lagrangian} subbundle; clearly $E$ is Lagrangian if and only if $E = E^{\perp}$.
\par
Conversely, one would like to know when a short exact sequence
\[ 0 \to E \to W \to E^{*} \to 0 \]
is induced by a symplectic form. We have
\begin{crit} An extension $0 \to E \to W \to E^{*} \to 0$ has a symplectic structure with respect to which $E$ is isotropic if and only if $W$ is isomorphic as a vector bundle to an extension whose cohomology class belongs to $H^{1}(X, \Sym^{2}E)$. \label{fundam}\end{crit}
\textbf{Proof}\\
Due to S.\ Ramanan; see \cite{Hit2005b}, $\S$ 2 for a proof.

\subsection*{Moduli of vector bundles over \textit{X}}

We recall the notions of stability and semistability for vector bundles. The \textsl{slope} of a bundle $W \to X$ is the ratio $\degree(W) / \rank(W)$, denoted $\mu(W)$. Then $W$ is \textsl{stable} (respectively, \textsl{semistable}) if for all proper nonzero subbundles $V \subset W$ we have
\[ \mu(V) < \mu(W) \quad \hbox{(respectively, $\mu(V) \leq \mu(W)$ ).} \]

Now let $r$ and $d$ be integers with $r \geq 1$. The moduli space of semistable vector bundles of rank $r$ and degree $d$ over $X$ is denoted $\urd$. For any line bundle $L \to X$ of degree $d$, the closed subvariety of $\urd$ of bundles with determinant $L$ is denoted $\surL$. References for these objects include Seshadri \cite{Ses1982} and Le Potier \cite{LeP1997}. The variety ${\cal U}_{X}(1,d)$ is the $d$th Jacobian variety of $X$, and will be denoted $J^{d}_X$. See for example Birkenhake--Lange \cite{BL2004}, Chap.\ 11, for details.
\par
The main object of interest for us is the moduli space of semistable symplectic vector bundles of rank $2n$ over $X$, which we denote $\Mn$. Since such a bundle always has trivial determinant, there is a forgetful map $\Mn \dashrightarrow {\cal SU}_{X}(2n, O_{X})$. We have

\begin{thm} $\Mn$ is canonically isomorphic to the moduli space of semistable principal $\Spn$-bundles over $X$, and the forgetful map $\Mn \to {\cal SU}_{X}(2n, O_{X})$ is an injective morphism. \label{Mn=Mxn}
\end{thm}
\textbf{Proof} (Sketch; see \cite{Hit2005a}, Chap.\ 1 for details.)\\
For any smooth variety $Y$, there is an equivalence between the groupoids
\[ \left\{ \hbox{principal $\Spn$-bundles over $Y$} \right\} + \{ \hbox{isomorphisms} \} \]
and
\begin{multline*} \left\{ \hbox{ symplectic vector bundles of rank $2n$ over $Y$ } \right\} \\ +  \{ \hbox{vector bundle isomorphisms respecting the symplectic structures} \}. \end{multline*}
Results of Ramanan \cite{Ram1981} and Ramanathan \cite{Rthn1975} on orthogonal bundles are easily adapted to the symplectic case to show that the notions of semistability and $S$-equivalence as principal $\Spn$-bundle and vector bundle coincide under this equivalence. This shows that the map $\Mn \to {\cal SU}_{X}(2n, O_{X})$ is an injective morphism. \qed \\
\\
\textbf{Remark:} In particular, any two symplectic forms on a polystable vector bundle differ by a vector bundle automorphism.\\
\\
Thm.\ \ref{Mn=Mxn} allows us to use results of Ramanathan \cite{Rthn1996a}, \cite{Rthn1996b} (especially Thm.\ 5.9) and \cite{Rthn1975} on moduli of principal $G$-bundles to deduce information about $\Mn$. We find that $\Mn$ is a projective variety of dimension $n(2n+1)(g-1)$ which is irreducible and normal. Moreover, the locus of stable vector bundles in $\Mn$ is dense.

\subsection*{Statement of the main theorem}

\textbf{Motivation:} Consider firstly the rank $2$ case. The following is well known, but we give a proof which illustrates Criterion \ref{fundam}.
\begin{prop} A bundle $W \to X$ of rank $2$ has an $O_X$-valued symplectic form if and only if it has trivial determinant. \end{prop}
\textbf{Proof}\\
If $W$ is symplectic of rank 2 then any line subbundle $L \subset W$ is clearly Lagrangian, so $W$ is an extension of $L^{-1}$ by $L$. Thus $W$ has trivial determinant.
\par
Conversely, if $\det ( W ) = O_X$ then any line subbundle $L \subset W$ induces an exact sequence $0 \to L \to W \to L^{-1} \to 0$, with class in $H^{1}(X, \Homom(L^{-1}, L))$. But since $L$ has rank $1$, we have $\Homom(L^{-1}, L) = \Sym^{2} L$, so $W$ is symplectic by Criterion \ref{fundam}.  \qed \\

This gives rise to an identification between ${\cal M}_1$ and ${\cal SU}_{X}(2, O_{X})$. When $X$ is nonhyperelliptic of genus $3$, Narasimhan and Ramanan \cite{NR1987} construct a generically finite cover of ${\cal M}_1$ (in this case isomorphic to the \textsl{Coble quartic}) by the union of the spaces $\PP H^{1}(X, L^{-2})$ as $L$ ranges over $J^{1}_X$. This description is a useful tool: for example, it was used by Pauly in \cite{Pau2002} to prove that the Coble quartic is self-dual.\\
\par
In the present paper we give a similar description of $\Mt$ when $X$ has genus $2$. The main tool is Criterion \ref{fundam}. We will consider a union of projective spaces of the form $\PP H^{1}(X, \Sym^{2}E)$ as $E$ ranges over a certain family of vector bundles over $X$.\\
\\
\textbf{Notation:} Let $V$ be a vector space and $v \in V$ a nonzero vector. Then we denote by $\langle v \rangle$ the point defined by $v$ in the projective space $\PP V$. We use similar notation for points of projectivised vector bundles.\\
\par
The moduli space $\utmo$ of semistable vector bundles of rank $2$ and degree $-1$ over $X$ is an irreducible variety of dimension $4g-3=5$. Since $\gcd(2,-1) = 1$, we have a diagram
\[ \xymatrix{ & \mathbf{E} \ar[d] & \\ & \utmo \times X \ar[dl]_{p} \ar[dr] & \\ \utmo & & X } \]
where $\mathbf{E}$ is a Poincar\'e bundle. By Riemann--Roch and semistability, the sheaf $R^{1}p_{*} ( \Sym^{2}\mathbf{E} )$ is locally free of rank $6$ over $\utmo$. We denote  $\mathbb{B}$ the associated $\PP^5$-bundle over $\utmo$. The fibre of $\mathbb{B}$ over a bundle $E$ is isomorphic to $\PP H^{1}(X, \Sym^{2}E)$; henceforth we denote this space $\PP^{5}_E$.
\par
By the results in Seshadri \cite{Ses1982}, Appendix II, there exists a vector bundle over $\mathbb{B} \times X$ whose restriction to $\{ ( \langle \delta \rangle ,E) \} \times X$ is isomorphic as a vector bundle to the extension $0 \to E \to W \to E^{*} \to 0$ defined by $\delta$. Thus there exists a classifying map $\Phi \colon \mathbb{B} \dashrightarrow \Mt$ by the moduli property of $\Mt$. A priori, $\Phi$ is only a rational map. The main result of this paper is
\begin{thm}
The moduli map $\Phi$ is a surjective morphism which is generically finite of degree $24$.
\label{main}
\end{thm}

The strategy for the proof is as follows. Firstly, we show that $\Phi$ is defined everywhere, so it is in fact a morphism. Then we consider a general (in a sense to be made precise) symplectic extension $0 \to E \to W \to E^{*} \to 0$ and show that it admits only finitely many isotropic subbundles of rank $2$ and degree $-1$. It is then easy to show that $\Phi$ is generically finite. By dimension count, irreducibility and projectivity of the spaces, it is surjective. Finally, we use results of Lange and Newstead \cite{LNew2002} to calculate the degree of $\Phi$.\\
\\
We conclude with some actual and future applications for this description of $\Mt$.

\paragraph{Acknowledgements:} I thank my doctoral supervisors C.\ Pauly, J.\ Bolton, W.\ Klingenberg, and W.\ Oxbury for their time and ideas. I am grateful to A.\ Beauville, A.\ Hirschowitz, K.\ Hulek (who pointed out Corollary \ref{KDMt}) and S.\ Ramanan (who showed me Criterion \ref{fundam}) for interesting and useful discussions. I acknowledge gratefully the financial and practical support of the Universities of Durham, Nice--Sophia-Antipolis and Hannover, the Engineering and Physical Sciences Research Council, and the Deutsche Forschungsgemeinschaft Schwerpunkt ``Globale Methoden in der Komplexen Geometrie''.

\section{Extensions and principal parts}

In this brief section we state some results on extensions of vector bundles which will be used throughout the paper. Firstly, we recall some notions from Kempf \cite{Kem1983}. The sheaf of sections of a vector bundle ($W$, $L$, $O_X$, $K_X$) will be denoted by the corresponding script letter (${\cal W}$, ${\cal L}$, $\strs$, $\cans$). For any vector bundle $V \to X$, there is a short exact sequence of $\strs$-modules
\[ 0 \to {\cal V} \to \sRat(V) \to \sPrin(V) \to 0 \]
where $\sRat(V)$ is the sheaf of rational sections of $V$ and $\sPrin(V)$ that of principal parts with values in $V$. Taking global sections, we get
\begin{equation} 0 \to H^{0}(X, V) \to \Rat(V) \to \Prin(V) \to H^{1}(X, V) \to 0. \label{cohomseq} \end{equation}
We write $\overline{s}$ for the principal part of $s \in \Rat(V)$, and we denote the cohomology class of $q \in \Prin(V)$ by $[q]$. Recall that an \textsl{elementary transformation} of a vector bundle $F \to X$ is the bundle corresponding to the kernel of a surjective map from ${\cal F}$ to a torsion sheaf.

\begin{thm}
Let $E$ and $F$ be vector bundles over $X$ and suppose that there are no nonzero maps $F \to E$. Let $0 \to E \to W \to F \to 0$ be an extension of class $\delta(W) \in H^{1}(X, \Homom(F,E))$.
\begin{enumerate}
\renewcommand{\labelenumi}{(\roman{enumi})}
\item There is a bijection between
\[ \left\{ \begin{array}{c} \hbox{principal parts } p \in \Prin(\Homom(F,E)) \\
\hbox{such that } \delta(W) = [p] \end{array} \right\} \]
and
\[ \left\{ \begin{array}{c} \hbox{elementary transformations of } $F$ \\
\hbox{lifting to vector subbundles of } W \end{array} \right\} \]
given by $p \leftrightarrow \Ker \left( p \colon {\cal F} \to \sPrin(E) \right)$.
\item If $F = E^*$ and $\delta(W)$ is symmetric, so that by Crit.\ \ref{fundam} we have a symplectic form on $W$, then the subbundle associated to $\Ker(p)$ is isotropic if and only if $p$ is itself a symmetric principal part\footnote{This is a stronger requirement than that the class $[p]$ be symmetric. By (\ref{cohomseq}), the class $[p]$ is symmetric if ${^{t}p} - p = \overline{\alpha}$ for some global rational section $\alpha$ of $\Homom(E^{*},E)$.}.
\end{enumerate}
\label{extprops} 
\end{thm}
\textbf{Proof}\\
See \cite{Hit2005b}, $\S$ 3.

\section{Vector bundles of rank 2 and odd degree over a curve of genus 2}

Here we assemble some technical results which will be needed later.

\subsection*{A cover of another moduli space}
The moduli space ${\cal U}_{X}(2, 2k+1)$ of semistable vector bundles of rank $2$ and degree $2k+1$ over $X$ is a variety of dimension $5$ which is smooth and irreducible.
\par
We give another description of ${\cal U}_{X}(2, 2k+1)$, which will in fact illustrate the strategy of the proof of our main result. In the same way as we arrived at the $\PP^5$-bundle $\mathbb{B} \to \utmo$, we construct a projective bundle $\mathbb{S} \to J^{k+1}_{X} \times J^{k}_X$ whose fibre at $(M,L)$ is isomorphic to $\PP H^{1}(X, \Homom(M,L)) = \PP^1$. This admits a moduli map $\Phi_{\mathbb{S}} \colon \mathbb{S} \dashrightarrow {\cal U}_{X}(2, 2k+1)$.
\begin{lemma}
The map $\Phi_{\mathbb{S}}$ is a surjective morphism of degree $4$. In particular, every $E \in {\cal U}_{X}(2, 2k+1)$ fits into $1$, $2$, $3$ or (generically) $4$ short exact sequences of the form
\begin{equation} 0 \to L \to E \to M \to 0 \label{extE} \end{equation}
where $L$ and $M$ are line bundles of degree $k$ and $k+1$ respectively.\label{lsubbsEth}
\end{lemma}
\textbf{Proof}\\
Let $E$ belong to the image of $\Phi_{\mathbb{S}}$, so $E$ is a stable extension of type (\ref{extE}). Suppose $L^{\prime} \subset E$ is a line subbundle of degree $k$. If $L^{\prime} \neq L$ then $L^{\prime}$ must be of the form $M(-x)$ for some point $x \in X$. Conversely, by Thm.\ \ref{extprops} (i) a subsheaf ${\cal M}(-x) \subset {\cal M}$ lifts to a line subbundle of $E$ if and only if $\delta(E)$ can be represented by a principal part which has a simple pole at $x$ and is otherwise regular.
\par
Now by analogy with Kempf--Schreyer \cite{KS1988}, $\S$ 1, we note that the standard map
\[ \phi_{|K_{X} M L^{-1}|} \colon X \dashrightarrow |K_{X} M L^{-1}|^{*} \cong \PP H^{1}(X, \Homom(M, L)) = \PP^{1} \]
can be defined by sending a point $x \in X$ to the line spanned by the cohomology class of a principal part which has a simple pole at $x$ and is otherwise regular, when this class is nonzero. Since $\phi_{|K_{X} M L^{-1}|}$ is surjective of degree at most 3 and $\delta(E)$ is nonzero, there are at most three, and generically three, points $x \in X$ such that $M(-x)$ is a subbundle of $E$.
\par
Therefore $E$ occurs as an extension of type (\ref{extE}) for at most four pairs $(M,L)$. Moreover, two distinct classes in the pencil $|K_{X} M L^{-1}|^*$ define nonisomorphic vector bundles: it is easy to check that there is at least one degree $k$ line subbundle which belongs to one but not to the other. Hence $\Phi_{\mathbb{S}}$ is quasifinite of degree $4$. By dimension count, it is dominant. Now it is not hard to see that every nontrivial extension of type (\ref{extE}) is stable, so $\Phi_{\mathbb{S}}$ is a morphism. Since $\mathbb{S}$ and ${\cal U}_{X}(2, 2k+1)$ are projective varieties, the image is closed, so it is surjective.
 \qed \\
\textbf{Note:} Lange and Narasimhan prove in \cite{LNar1983}, Prop.\ 4.2, that a vector bundle of rank $2$ and degree $2k+1$ over $X$ has at most $4$ subbundles of degree $k$.\\
\\
We will need
\begin{lemma}
Let $V$ and $L$ be vector bundles of ranks $n$ and $1$ respectively, and let $f \colon L \to V$ be a homomorphism. Then $f$ factorises via a map $L(x) \to V$ if and only if $f$ is zero at $x$. \label{injectivity} \end{lemma}
\textbf{Proof}\\
Narasimhan--Ramanan \cite{NR1969}, Lemma 5.3.

\begin{prop}
Let $E \to X$ be a (semi)stable bundle of rank $2$ and degree $2k+1$. Then $h^{0}(X, \Homom(L, E)) \leq 1$ for all $L \in J^{k}_X$.
\label{oneinj} \end{prop}
\textbf{Proof}\\
Let $E$ be a vector bundle of rank $2$ and degree $2k+1$ and suppose that $h^{0}(X, \Homom(L, E)) \geq 2$ for some $L \in J^{k}_X$. Firstly, note that if any nonzero map $L \to E$ has a zero then $E$ would have a line subbundle of degree at least $k+1$ by Lemma \ref{injectivity}, so would not be semistable. Thus we can suppose that $E$ is an extension $0 \to L \to E \to M \to 0$ for some $M \in J^{k+1}_X$. By hypothesis, there also exists a copy of $L$ in $E$ whose projection to $M$ is generically nonzero, so $L = M(-x)$ for some $x \in X$. Therefore, by Thm.\ \ref{extprops} (i) the class $\delta(W)$ can be represented by a principal part $p$ with just one pole which is simple and at $x$.
\par
But since $L = M(-x)$, the bundle $M^{-1}L(x)$ is trivial, so has a global section which is nonzero at $x$. This gives a global rational section of $M^{-1}L$ with just a simple pole at $x$. Hence $[p] = 0$ by (\ref{cohomseq}), so $E$ is a trivial extension. In particular, it is not semistable. \qed

\subsection*{Genericity of bundles in $\utmo$}

In this section, we show that some conditions which we will need later are satisfied by a generic stable bundle of degree $-1$ and rank $2$ over $X$. We will need
\begin{lemma}
Every semistable vector bundle $F$ of rank at most $3$ and slope $1$ over a curve of genus $2$ satisfies $h^{0}(X, M \otimes F) = 0$ for generic $M \in J^{0}_X$. \label{g2rk3thetadiv}
\end{lemma}
\textbf{Proof}\\
This is a special case of Raynaud \cite{Ray1982}, Cor.\ 1.7.4.
 
\begin{prop}
For generic $E \in \utmo$, the bundle $\Homom(E, E^{*})$ has no global sections. \label{etegs}
\end{prop}
\textbf{Proof}\\ We have $h^{0}(X, \Homom(E, E^{*})) = h^{0}(X, \det E^{*}) + h^{0}(X, \Sym^{2}E^{*})$, so it suffices to show that the subsets of $\utmo$ where $h^{0}(X, \det E^{*}) > 0$ and where $h^{0}(X, \Sym^{2}E^{*}) > 0$ are each of codimension $1$.
\par
Firstly, $h^{0}(X, \det E^{*}) > 0$ if and only if $(\det E)^{-1}$ is effective. The set of such $E$ is the inverse image of $\Supp(\Theta)$ under the map $\utmo \to J^{1}_X$ defined by $E \mapsto (\det E)^{-1}$. Clearly this map is surjective, so the inverse image of a divisor is a divisor.
\par
For the rest: $\Sym^{2}E^*$ is of slope $1 = g-1$ for all $E \in \utmo$, so we expect that if the set
\[ \left\{ E \in \utmo : h^{0}(X, \Sym^{2} E^{*}) > 0 \right\} \]
is not equal to $\utmo$ then it is the support of a divisor. Since $\utmo$ is irreducible, then, it suffices to exhibit one $E$ such that $h^{0}(X, \Sym^{2}E^{*}) = 0$. Choose any $E \in \utmo$. If $h^{0}(X, \Sym^{2}E^{*}) = 0$ then we are done. If not, we note that $\Sym^{2}E$ is semistable for all $E \in \utmo$ by Le Potier \cite{LeP1997}, p.\ 161. Thus, by Lemma \ref{g2rk3thetadiv}, there exists at least one $M \in J^{0}_X$ such that
\[ h^{0}(X, M \otimes \Sym^{2}E^{*}) = 0. \]
Let $N$ be any square root of $M$. Since
\[ M \otimes \Sym^{2}E^{*} = N^{2} \otimes \Sym^{2}E^{*} \cong \Sym^{2}(N \otimes E^{*}), \]
by construction the bundle $E^{\prime} := N^{-1} \otimes E$, which is stable of degree $-1$ and rank $2$, satisfies $h^{0}(X, \Sym^{2}(E^{\prime})^{*}) = 0$. This completes the proof of the proposition.
 \qed \\
\\
Recall that the \textsl{Weierstrass points} of $X$ are the 6 points fixed by the hyperelliptic involution $\iota \colon X \to X$.

\begin{prop} A general $E \in \utmo$ satisfies the following:
\begin{itemize}
\item $E^*$ has no degree $1$ line bundle quotients $L$ such that $L^{2} = O_{X}(2x)$. \label{genEL}
\item $\det ( E^{*} )^2$ is of the form $O_{X}(a+b)$ for distinct $a, b \in X$.
\item Furthermore, neither $a$ nor $b$ is a Weierstrass point.
\end{itemize} \label{genE} \end{prop}
\textbf{Proof}\\
We define some hypersurfaces in $J^{1}_{X} \times J^{0}_X$. Firstly, let $Y$ be the set
\[ \left\{ L \in J^{1}_{X} : \hbox{ $L^2$ is of form $O_{X}(2x)$ for some $x \in X$} \right\}, \]
the union of all translates of the theta divisor by line bundles of order 2 in $J^{0}_X$. We let $H_{1} := Y \times J^{0}_X$.
\par
Next, consider the curve $X_{2} := \{ O_{X}(2x) : x \in X \}$ in $J^{2}_X$. Let $H_2$ be the inverse image in $J^{1}_{X} \times J^{0}_X$ of $X_2$ by the map $J^{1}_{X} \times J^{0}_{X} \to J^{2}_X$ given by $(L, M) \mapsto L^{2}M^2$.
\par
Finally, for each Weierstrass point $w \in X$ we have the set
\[ \{ O_{X}(w + x) : x \in X \} =: \Theta + w. \]
Write $H_3$ for the inverse image of the union of the six $\Theta + w$ by the map $(L, M) \mapsto L^{2} M^2$.
\par
Let $H := H_{1} \cup H_{2} \cup H_3$. Now recall from Lemma \ref{lsubbsEth} that we have a cover $\Phi_{\mathbb{S}} \colon \mathbb{S} \to \uto$ (here $k=0$). If any of the above conditions were not satisfied by a general $E \in \utmo$ then the restriction of $\Phi_{\mathbb{S}}$ to the subvariety $\mathbb{S}|_{H} \subset \mathbb{S}$ would be dominant. But this is impossible since $\dim ( \mathbb{S}|_{H} ) = 4$  but $\dim(\uto) = 5$. The proposition follows.  \qed

\subsection*{A ruled surface in $\PP^5$}

In this section, we fix $E \in \utmo$ and describe a map of the ruled surface $\PP E$ into $\PP^{5}_E$. We then study the extent to which this map fails to be an embedding for generic $E$.
\par
We generalise slightly the aforementioned approach from $\S$ 1 of Kempf--Schreyer \cite{KS1988}. For any $x \in X$, we can form the exact sheaf sequence
\[ 0 \to \Sym^{2}{\cal E} \to \left( \Sym^{2}{\cal E} \right) (x) \to \frac{\left( \Sym^{2}{\cal E} \right) (x)}{\Sym^{2}{\cal E}} \to 0 \]
whose cohomology sequence begins
\[ 0 \to H^{0}(X, (\Sym^{2}E)(x)) \to \left( \Sym^{2}E \right) (x)|_{x} \xrightarrow{\delta} H^{1}(X, \Sym^{2}E) \to \ldots \]
The second term can be identified with the set of $\Sym^{2}E$-valued principal parts which are regular except for possibly a simple pole at $x$. Moreover, there is a canonical isomorphism $\PP ( \Sym^{2}E )(x)|_{x} \xrightarrow{\sim} \PP \: \Sym^{2}E|_x$.
\par
For each $x \in X$, then, we define a map $\psi_{x} \colon \PP E|_{x} \to \PP H^{1}(X, \Sym^{2}E)$ by the composition
\begin{equation} \PP E|_{x} \xrightarrow{\mathrm{Segre}} \PP \: \Sym^{2}E|_{x} \xrightarrow{\sim} \PP ( \Sym^{2}E ) (x)|_{x} \xrightarrow{\PP \delta} \PP H^{1}(X, \Sym^{2} E). \label{defpsi}\end{equation}
We define $\psi \colon \PP E \to \PP H^{1}(X, \Sym^{2}E)$ to be the product of the $\psi_x$ over all $x \in X$. We will need
\begin{prop} For general $E \in \utmo$, the bundle $K_{X} \otimes \Sym^{2}E^{*}$ is generated by global sections. \label{genEglobsect} \end{prop}
\textbf{Proof}\\
For any $x \in X$, we have an exact cohomology sequence
\begin{multline*} 0 \to H^{0}(X, K_{X}(-x) \otimes \Sym^{2}E^{*}) \to H^{0}(X, K_{X} \otimes \Sym^{2}E^{*}) \to \\
\left( K_{X} \otimes \Sym^{2}E^{*} \right)|_{x} \to H^{1}(X, K_{X}(-x) \otimes \Sym^{2}E^{*}) \to 0, \end{multline*}
whence $K_{X} \otimes \Sym^{2}E^{*}$ is generated by global sections if and only if
\[ h^{1}(X, K_{X}(-x) \otimes \Sym^{2}E^{*}) = h^{1}(X, (\Sym^{2}E^{*})(\iota (x))) = 0 \]
for all $x \in X$. By Serre duality, this is equal to
\begin{align*} h^{0}(X, K_{X}(- \iota(x)) \otimes \Sym^{2}E) &= h^{0}(X, (\Sym^{2}E)(x)) \\
&= h^{0}(X, \Homom ( O_{X}(-x) , \Sym^{2}E)). \end{align*}

Now let $M \in J^{0}_X$ be a bundle which is not of order 2, and $N \in J^{-1}_X$ such that $N^{-1}M^{-1}$ is not effective; clearly the general $(N, M) \in J^{-1}_{X} \times J^{0}_X$ satisfies these conditions. We will show that for all but at most two extensions
\[ 0 \to N \to E \xrightarrow{c} M \to 0, \]
there are no maps $O_{X}(-x) \to \Sym^{2} E$ for any $x \in X$. Let $E$ be such an extension, with class $p_{1} + p_{2} + p_{3} \in |K_{X}MN^{-1}|$. It is not hard to see that the map $E \otimes E \to M \otimes E$ given by $e \otimes f \mapsto \frac{1}{2} ( c(e) \otimes f + c(f) \otimes e)$ induces another exact sequence
\[ 0 \to N^{2} \to \Sym^{2} E \to M \otimes E \to 0. \]
Since $M$ is not of order 2, there is a unique pair of points $q_{1}, q_{2} \in X$ such that $K_{X}M^{2} = O_{X}(q_{1} + q_{2})$.
\par
If there is a nonzero map $O_{X}(-x) \to \Sym^{2}E$ then $O_{X}(-x)$ must be a subbundle of $M \otimes E$, equivalently $M^{-1}(-x)$ must be a subbundle of $E$. But the degree $-1$ subbundles of $E$ are exactly $N$ and the $M(-p_{i})$. If $N=M^{-1}(-x)$ then $N^{-1} M^{-1}$ is effective, contrary to hypothesis.
\par
If $M(-p_{i})=M^{-1}(-x)$ for some $i$ then $M^{2} = O_{X}(p_{i} - x)$, so we obtain $K_{X}M^{2} = O_{X}(p_{i} + \iota(x))$. This means that the only extensions of $M$ by $N$ such that $M^{-1}(-x)$ belongs to $E$ are those where $p_{i} = q_j$ for some $i$ and $j$. This can happen for at most two classes in $|K_{X}MN^{-1}|$.
\par
In summary, for all but at most two extensions of $M$ by $N$, there are no maps $O_{X}(-x) \to \Sym^{2}E$. The proposition follows.  \qed \\

Now write $\pi$ for the projection $\PP E \to X$ and let $\Upsilon \to \PP E$ be the line bundle $\pi^{*}K_{X} \otimes {\cal O}_{\PP E}(2)$.
\begin{prop} There is a natural identification $\PP^{5}_{E} \cong |\Upsilon|^{*}$ under which $\psi$ coincides with the natural map $\phi_{|\Upsilon|} \colon \PP E \dashrightarrow |\Upsilon|^*$. In particular, $\psi$ is algebraic. Moreover, for general $E$, it is defined everywhere. \label{phialg}
\end{prop}
\textbf{Proof} (Sketch; see \cite{Hit2005a}, Prop.\ 4.11 for details)\\
Firstly, using the fact that $K_X$ and $\Sym^{2}E^{*}$ are locally trivial on $X$, one sees that sections of $\Upsilon$ over $\PP E$ can be interpreted naturally as sections of $K_{X} \otimes \Sym^{2}E^{*}$ over $X$. Thus for any $u \in \PP E$, a principal part $p$ defining $\psi(u)$ defines a linear form $[p]$ on $H^{0}(\PP E, \Upsilon)$ by Serre duality. One shows that $\Ker [p]$ contains the subspace $H^{0}(\PP E, \Upsilon - u)$ of sections vanishing at $u$. Now by Prop.\ \ref{genEglobsect}, for general $E$, the bundle $K_{X} \otimes \Sym^{2}E^*$ is generated by global sections, so there is a section of $\Upsilon$ not vanishing at $u$. By this and the fact that no global rational section of $K_X$ has just one simple pole, we find that $[p]$ is nonzero. Hence $H^{0}(\PP E, \Upsilon - u)$ is a hyperplane. It follows that $\psi$ is identified with $\phi_{|\Upsilon|}$ and is a morphism.
 \qed \\
\\
We now determine the extent to which $\psi$ fails to be an embedding in general.
\begin{lemma}
For general $E \in \utmo$, the map $\psi$ is an embedding except at a finite number of points. \label{psiemb} \end{lemma}
\textbf{Proof}\\ Since $E$ is general, we can suppose that $\psi$ is a morphism by Props.\ \ref{genEglobsect} and \ref{phialg}. We distinguish three ways in which it can fail to be an embedding:
\begin{enumerate}
\renewcommand{\labelenumi}{(\roman{enumi})}
\item $\psi(u) = \psi(v)$ for some $u, v \in \PP E$ lying over distinct $x, y \in X$.
\item $\psi(u) = \psi(v)$ for distinct $u$ and $v$ in a fibre $\PP E|_x$.
\item The differential of $\psi$ is not injective at a point $u$ in some $\PP E|_x$.
\end{enumerate}
Recall that we have the cohomology sequence
\begin{equation} 0 \to \Rat(\Sym^{2}E) \to \Prin(\Sym^{2}E) \xrightarrow{\delta} H^{1}(X, \Sym^{2}E) \to 0 \label{S2Ecohom} \end{equation}
which is exact since $h^{0}(X, \Sym^{2}E) = 0$ by semistability.
\par
Suppose that (i) occurs. Let $\tilde{u}$ and $\tilde{v}$ be lifts of $u$ and $v$ to $E$ and suppose that $p$ and $q$ are $\Sym^{2}E$-valued principal parts with simple poles along $\tilde{u} \otimes \tilde{u}$ and $\tilde{v} \otimes \tilde{v}$ respectively. Then
\[ \langle [ p ] \rangle = \psi(u) = \psi(v) = \langle [ q ] \rangle. \]
After normalising if necessary, we can assume that the class $[ p - q ]$ is zero in $H^{1}(X, \Sym^{2}E)$. Then by exactness of (\ref{S2Ecohom}), there exists a global rational section $\alpha$ of $\Sym^{2}E$ with principal part $p - q$. This is a global regular section of $(\Sym^{2}E)(x+y)$ with $\alpha(x)$ and $\alpha(y)$ decomposable. The generic rank of the corresponding symmetric map $\alpha \colon E^{*} \to E(x+y)$ may be 1 or 2. If it is 2 then $\det ( \alpha )$ is generically nonzero, so has
\[ \degree ( E^{*} ) - \degree \left( E(x+y) \right) = 2 \]
zeroes. It follows that $\det(\alpha) \in H^{0}(X, \strs(x+y))$. But we also have
\[ \det(\alpha) \in H^{0}(X, \Homom( \det E^{*}, (\det E)(2x+2y))), \]
so $(\det E)^{2} = \strs(-x-y)$.
\par
Now since $E^*$ is of rank $2$, there is an isomorphism $E \xrightarrow{\sim} E^{*} \otimes \det ( E )$, which induces another isomorphism $\Sym^{2}E \xrightarrow{\sim} \Sym^{2}E^{*} \otimes \det ( E )^2$. Thus
\[ \left( \Sym^{2}E \right) (x+y) \cong \Sym^{2}E \otimes (\det E)^{-2} \cong \Sym^{2}E^{*}, \]
so $\alpha$ defines a nonzero section of $\Sym^{2}E^*$. But by Prop.\ \ref{etegs}, for general $E$ there are no such $\alpha$.
\par
Therefore we can suppose $\alpha$ to be generically of rank 1. This means that $\alpha$ factorises as
\[ {\cal E}^{*} \to {\cal L} \to {\cal E}(x+y) \]
where ${\cal L}$ is an invertible subsheaf of ${\cal E}(x+y)$. By stability of $E^*$, the degree of ${\cal L}$ is at least 1, and by stability of $E(x+y)$, it is at most 1. Therefore it is 1 and in fact we have a sequence of vector bundle maps $E^{*} \to L \to E(x+y)$. Since $\alpha$ is symmetric, this is identified with
\[ E^{*} \to L^{-1}(x+y) \to E(x+y), \]
so $L^{2} = O_{X}(x+y)$. We see that $L^{-1}$ is a subbundle of $E$ and the two points of $\PP E$ that are contracted are those corresponding to the fibres $L^{-1}|_x$ and $L^{-1}|_y$.
\par
In summary, for generic $E$, the map $\psi$ contracts two points of $\PP E$ only if they are the images of $N|_x$ and $N|_y$ for some degree $-1$ subbundle $N$ of $E$ such that $N^{-2} = O_{X}(x+y)$. By Lemma \ref{lsubbsEth}, there are at most four such $N$. By Prop.\ 4, we can assume that $N^{-2} \neq K_X$ (otherwise the quotient $N^{-1}$ of $E^*$ would satisfy $(N^{-1})^{2} = O_{X}(2w)$ for any Weierstrass point $w$), so there is a unique pair $x,y$ for each $N$. Moreover, by Prop.\ \ref{oneinj}, there is only one independent map $N \to E$, so there is only one possibility for the points $N|_x$ and $N|_y$ in $\PP E$.
\par
Putting all this together, we have at most eight points of $\PP E$ at which $\psi$ fails to be injective in this way.\\
\par
Case (ii) is simpler. Proceeding as above, we get a symmetric map $\alpha \colon E^{*} \to E(x)$ which is of rank 2 at $x$. Since the bundles have the same rank and degree, it is an isomorphism. Twisting $\alpha^{-1}$ by $O_{X}(-x)$, we get a symmetric map $E \to E^{*}(-x)$. Composing with the generic inclusion $E^{*}(-x) \to E^{*}$, we obtain a nonzero global section of $\Sym^{2}E^*$ so, by Prop.\ \ref{etegs}, the bundle $E$ is not general.\\
\par
Lastly, suppose (iii) happens. This means that the space of global sections of $\Upsilon$ which vanish to order 2 at $u$ is of dimension greater than expected, that is, at least 4. Let $e$, $f$ be a local basis for $E$ with $\langle e|_{x} \rangle = u$, and consider also the dual local basis $e^{*}, f^{*}$ of $E^*$. Let $z$ be a local coordinate on $X$ centred at $x$. Near $u$, a section of $\Upsilon$ vanishing to order 2 at $u$, viewed as a section of $K_{X} \otimes \Sym^{2} E^*$, has the form
\begin{multline*} dz \otimes \left( c_{0} ( f^{*} \otimes f^{*} ) + b_{1} z ( e^{*} \otimes f^{*} + f^{*} \otimes e^{*} ) + c_{1} z ( f^{*} \otimes f^{*} ) \right. \\
\left. + \hbox{ terms of higher order in $z$} \right) \end{multline*}
for some $c_{0},  b_{1}, c_{1} \in \C$. Now clearly the contractions of such a section against the principal parts
\[ \frac{e \otimes e}{z} , \quad \frac{e \otimes f + f \otimes e }{z} \quad \hbox{and} \quad \frac{e \otimes e}{z^2} \]
are all regular, so the cohomology classes of these principal parts belong to the space\footnote{The projectivisation of this space is the embedded tangent space to $\psi ( \PP E )$ at $\psi(u)$.} $\Ker \left( H^{0}( \PP E , \Upsilon )^{*} \to H^{0}( \PP E , \Upsilon - 2u )^{*} \right)$. By hypothesis, this has dimension at most 2, so, by exactness of (\ref{S2Ecohom}), there is a relation
\[ \lambda_{1} \left[ \frac{e \otimes e}{z} \right] + \lambda_{2} \left[ \frac{e \otimes f + f \otimes e }{z} \right] + \lambda_{3} \left[ \frac{e \otimes e}{z^2} \right] = 0 \]
for some $\lambda_{1}, \lambda_{2}, \lambda_{3} \in \C$, not all zero. This means that there is a global rational section of $\Sym^{2}E$ with principal part
\[ \frac{ \lambda_{1}e \otimes e}{z} + \frac{ \lambda_{2} ( e \otimes f + f \otimes e ) }{z} + \frac{\lambda_{3}e \otimes e}{z^2}. \]
If $\lambda_{3} = 0$ then in fact we are in situation (ii), so $E$ is not general by Prop.\ \ref{etegs}. Assuming the contrary, there is a regular symmetric map $\alpha \colon E^{*} \to E(2x)$ which has rank 1 at $x$. Again, there are two possibilities: the generic rank could be either 1 or 2. If it is 1 then as before $\alpha$ factorises as $E^{*} \to L \to E(2x)$ where $L$ is a degree 1 line bundle quotient of $E^*$. Since $\alpha$ is symmetric, we have as before $L^{2} = O_{X}(2x)$, so $E$ is not general, this time by Prop.\ \ref{genE}.
\par
Thus we can suppose that $\alpha$ has generic rank 2. Then we have the sheaf sequence
\[ 0 \to {\cal E}^{*} \xrightarrow{\alpha} {\cal E}(2x) \to T \to 0 \]
where $T$ is a torsion sheaf. Now the determinant of $\alpha$ vanishes at two points (counted with multiplicity) because $\degree(E(2x)) - \degree ( E^{*} ) = 2$. One of these points is $x$; call the other one $y$ (they coincide if and only if $\lambda_{2} = 0$). So $\det ( \alpha ) \in H^{0}(X, O_{X}(x+y))$. On the other hand,
\[ \det ( \alpha ) \in H^{0}(X, \Homom( \det E^{*}, \det (E(2x)) )) = H^{0}(X, (\det E)^{2}(4x)). \]
By genericity and Prop.\ \ref{genE}, we have $\det (E^{*})^{2} = O_{X}(a+b)$ for distinct $a, b \in X$, neither one a Weierstrass point. Then $(\det E)^{2}(4x) = O_{X}(4x-a-b) = O_{X}(x+y)$, so $O_{X}(3x) = O_{X}(y+a+b)$.
\par
We show that there are only a finite number of solutions $x, y \in X$ to this equation. Otherwise, it is not hard to see that for each $x \in X$ we would have $y(x) \in X$ such that
\[ O_{X}(3x) = O_{X}(y(x)+a+b) \]
($y(x)$ is well defined if it exists because if $y+a+b$ and $y^{\prime}+a+b$ are linearly equivalent then clearly $y = y^{\prime}$). In particular, $O_{X}(3a) = O_{X}(y(a)+a+b)$. Then either $b+y(a) = 2a$, contradicting the hypothesis $b \neq a$, or $a$ is a base point of $O_{X}(3a)$, whence $a$ is a Weierstrass point, again contradicting genericity. Thus there are only finitely many $x \in X$ such that $y$ and $\alpha$ with these properties can exist.
\par
We show that if $y$ does exist, there is at most one independent $\alpha$. For otherwise,
\[ \det \colon H^{0}(X, \Sym^{2}E(2x)) \to H^{0}(X, O_{X}(x + y)) \]
has positive-dimensional fibres, so there is a nonzero symmetric homomorphism $\beta \colon E^{*} \to E(2x)$ with determinant everywhere zero, that is, of rank generically 1. But we saw above that, by genericity, this is impossible.\\
\par
In summary, for general $E$, the differential of $\psi$ can fail to be injective at only finitely many points of $\PP E|_x$. This completes the proof of the lemma.
 \qed \\

\section{Proof of the main theorem}

\noindent In this section we will prove Thm.\ \ref{main}. We begin by studying some special loci of $\PP^{5}_E$.

\subsection*{Stability of bundles in $\PP^{5}_E$}

We investigate the loci of strictly semistable bundles in $\PP^{5}_E$ and prove that the general extension represented in $\mathbb{B}$ is a stable vector bundle. We begin by stating a technical result:

\begin{prop}
Let $W \to X$ be any self-dual vector bundle. Then $W$ is stable (resp., semistable) if and only if it contains no destabilising (resp., desemistabilising) subbundles of rank at most $\frac{1}{2} \rank (W)$. \label{selfdualstable}
\end{prop}
\textbf{Proof}\\
This is straightforward to check; we remark that it is true whether $W$ has even or odd rank.

\begin{lemma}
For any $E \in \utmo$, every nontrivial symplectic extension $W$ of $E^*$ by $E$ is semistable. \label{evdef}
\end{lemma}
\textbf{Proof}\\
Since a symplectic vector bundle is self-dual, by Prop.\ \ref{selfdualstable} it is enough to show that $W$ contains no desemistabilising subbundles of rank at most \nolinebreak $2$.
\par
Firstly, suppose $L \subset W$ is a line subbundle of degree at least $1$. Since $E$ is stable, $h^{0}(X, \Homom(L, E^{*})) > 0$. But since $L$ and $E^*$ are stable and $\mu(L) \geq 1 > \mu(E^{*}) = \frac{1}{2}$, this is impossible.
\par
Next, suppose $F \subset W$ is a subbundle of rank $2$ and degree at least $1$. Since we have just seen that $W$ contains no line subbundles of positive degree, $F$ must be stable. Firstly, if the composed map $F \to W \to E^*$ were zero then we would have a nonzero map $F \to E$, contradicting stability of $F$ and $E$. Therefore the composed map is nonzero. Since $F$ and $E^*$ are stable, it must be an isomorphism, so $W$ is a trivial extension.
 \qed \\
\\
This lemma shows in particular that the classifying map $\Phi \colon \mathbb{B} \to \Mt$ is defined everywhere, so it is a morphism.\\
\\
We quote a result of Narasimhan and Ramanan:
\begin{lemma} \label{extlift}
Consider a diagram of vector bundles over $X$
\[ \xymatrix{ 0 \ar[r] & E \ar[r] & W \ar[r] & F \ar[r] & 0 \\ & & & V \ar[u]_{f} & } \]
where the top row is exact. Then $f$ factorises via a map $V \to W$ if and only if the class $\delta(W)$ of the extension belongs to the kernel of the induced map
\[ f^{*} \colon H^{1}(X, \Homom(F,E)) \to H^{1}(X, \Homom(V, E)). \]
\end{lemma}
\textbf{Proof}\\
Narasimhan--Ramanan \cite{NR1969}, Lemma 3.1.

\begin{lemma}
The generic symplectic bundle in $\mathbb{B}$ is a stable vector bundle. \label{genWstable}
\end{lemma}
\textbf{Proof}\\
We will prove this by determining the classes in a general $\PP^{5}_E$ representing vector bundles which contain a subbundle of degree $0$; this will also be used later. Again, by Prop.\ \ref{selfdualstable} it suffices to check for degree $0$ subbundles of ranks $1$ or $2$.
\par
Firstly, suppose $M \subset W$ is a line subbundle of degree $0$. Since $E$ is stable, $h^{0}(X, \Homom(M, E^{*})) > 0$; then in fact $M$ is a subbundle of $E^*$ by Lemma \ref{injectivity} since $E^*$ is stable. By Lemma \ref{lsubbsEth}, there are at most $4$ possibilities for $M \in J^{0}_X$.
\par
Conversely, by Lemma \ref{extlift}, a map $j \colon M \hookrightarrow E^*$ lifts to $W$ if and only if $\delta(W)$ belongs to the kernel of the induced map
\[ j^{*} \colon H^{1}(X, \Sym^{2}E) \to H^{1}(X, \Homom(M, E)). \]
We check that this map is surjective. By Serre duality, it is equivalent to check that the transposed map $H^{0}(X, K_{X} M \otimes E^{*}) \to H^{0}(X, K_{X} \otimes \Sym^{2}E^{*})$, which by abuse of notation we also denote $j$, is injective. Now the induced map
\[ j \colon H^{0}(X, K_{X} M \otimes E^{*}) \to H^{0}(X, K_{X} \otimes E^{*} \otimes E^{*}) \]
is injective, because the functor $\cans \otimes - \otimes {\cal E}^*$ and the global section functor are left exact (the latter because these sheaves are locally free). Since $H^{0}(X, K_{X} \otimes E^{*} \otimes E^{*})$ is the direct sum
\[ H^{0}(X, K_{X} \otimes \Sym^{2}E^{*}) \oplus H^{0}(X, K_{X} \otimes \bigwedge^{2} E^{*}), \]
we have to show that $\Image(j) \cap H^{0}(X, K_{X} \otimes \bigwedge^{2}E^{*}) = \{ 0 \}$. Now $E$ is of rank $2$, so the latter space is just $H^{0}(X, K_{X} \otimes (\det E)^{-1})$, and any nonzero section therein has $3$ zeroes (counted with multiplicity). But if a map $K_{X}^{-1}M^{-1} \to E^*$ vanished at $3$ points, $E^*$ would contain a line subbundle of degree at least $1$ by Lemma \ref{injectivity}; this would contradict the stability of $E^*$.
\par
Thus the restriction of $j^*$ to $H^{1}(X, \Sym^{2}E)$ is surjective. Since $\Sym^{2}E$ and $\Homom(M, E)$ have no global sections, the kernel of $j^*$ is of dimension
\[ -\chi(X, \Sym^{2}E) + \chi( \Homom(M, E)) = 6 - 3 = 3. \]
Hence there is a union of between $1$ and $4$ projective planes in $\PP^{5}_E$ representing extensions which are destabilised by a line subbundle of degree $0$.\\
\par
Now we consider a destabilising subbundle $G \subset W$ of rank $2$. We will use the map $\psi \colon \PP E \to \PP^{5}_E$ defined in section 3.
\begin{prop}
Let $E \in \utmo$ be general, and let $W$ be a nontrivial symplectic extension of $E^*$ by $E$. Then $W$ is destabilised by a subbundle of rank $2$ and degree $0$ if and only if $\langle \delta(W) \rangle$ belongs to $\psi(\PP E)$. \label{strsemist} \end{prop}
\textbf{Proof}\\
Let $G \subset W$ be a subbundle of rank $2$ and degree $0$. Then the rank of ${\cal G} \cap {\cal E}$ is $2$, $1$ or $0$. It cannot be $0$ because then the image of ${\cal G}$ in ${\cal E}^*$ would be a torsion subsheaf of length $1$, of which there are none. If ${\cal G} \cap {\cal E}$ were an invertible subsheaf ${\cal L}$ then we would have a diagram
\[ \xymatrix{ 0 \ar[r] & {\cal E} \ar[r] & {\cal W} \ar[r] & {\cal E}^{*} \ar[r] & 0 \\ 0 \ar[r] & {\cal L} \ar[r] \ar @{^{(}->}[u] & {\cal G} \ar[r] \ar @{^{(}->}[u] & {\cal M} \ar @{^{(}->}[u] \ar[r] & 0 } \]
where ${\cal M}$ is an invertible subsheaf of ${\cal E}^*$. Since $E$ and $E^*$ are stable, we have $\degree ( {\cal L} ) \leq -1$ and $\degree ( {\cal M} ) \leq 0$. But then $\degree ( G ) \leq -1$, a contradiction.
\par
Hence $G$ is an elementary transformation $0 \to {\cal G} \to {\cal E}^{*} \to \C_{x} \to 0$ for some point $x \in X$. Therefore, by Thm.\ \ref{extprops} (i), the class $\delta(W)$ is of the form $[q]$ where $q \in \Prin(\Homom(E^{*}, E)) \cong \Homom \left( {\cal E}^{*}, \sPrin(E) \right)$ has kernel ${\cal G}$. Clearly $q$ is supported at $x$ with a simple pole along $f \otimes f^{\prime}$ for some nonzero $f, f^{\prime} \in E|_x$.
\par
Now since $W$ is symplectic, ${^{t}q} - q = \overline{\alpha}$ for some global rational section $\alpha$ of $\Homom(E^{*}, E)$. Since $\overline{\alpha}$ is clearly antisymmetric, $\frac{\alpha - {^{t}\alpha}}{2}$ defines a global rational section of $\det E$ with principal part $\overline{\alpha}$. Now if $\alpha$ is nonzero then it has a single simple pole at $x$, so is a global regular section of $(\det E)(x)$. But this is nonzero only if $\det E = \strs(-x)$, which it is not since $E$ is general. Therefore $f^{\prime}$ is proportional to $f$ and $\langle \delta(W) \rangle = \psi \langle f \rangle$.\\
\par
Conversely, suppose that $\delta(W)$ lies over $\psi \langle f \rangle$ for some nonzero $f \in E|_x$. By the definition of $\psi$, the class $\delta(W)$ can be represented by $q \in \Prin(\Sym^{2}E)$ supported at one point $x \in X$ and with a simple pole along $f \otimes f$. By Theorem \ref{extprops} (i), the kernel of $q$ lifts to a subbundle of $W$ which clearly has rank $2$ and degree $0$.
 \qed \\
\\
In summary, we have shown that the locus of extensions in a general $\PP^{5}_E$ containing a subbundle of degree $0$ is of dimension $2$. The complement of this locus in $\PP^{5}_E$ consists of classes defining stable vector bundles. Lemma \ref{genWstable} follows.
 \qed

\subsection*{Maximal Lagrangian subbundles}

Let $E \in \utmo$ be general and consider a general stable symplectic extension $W$ of $E^*$ by $E$. We will show that $W$ has finitely many Lagrangian subbundles of degree $-1$. Suppose $F \subset W$ is such a subbundle. There are three possibilities for the rank of the sheaf ${\cal F} \cap {\cal E}$: these are $2$, $1$ and $0$. If it is $2$ then it is not hard to see that $F = E$.

\begin{prop}
A general extension in $\PP^{5}_E$ contains no subbundles of rank $2$ and degree $-1$ intersecting $E$ generically in rank $1$.
\label{FintEdrk1}
\end{prop}
\textbf{Proof} (idea from Lange--Newstead \cite{LNew2002}, Prop.\ 2.4)\\
Let $F \subset W$ be such a subbundle; if $W$ is stable then clearly $F$ is too. We write ${\cal L}$ for the subsheaf ${\cal F} \cap {\cal E}$, which is invertible because, say, ${\cal E}$ is locally free. By stability, $\degree ( {\cal L} ) \leq -1$. The image ${\cal M}$ of ${\cal F}$ in ${\cal E}^*$ is coherent, hence invertible because ${\cal E}^*$ is locally free. Since $E^*$ is stable, $\degree ( {\cal M} ) \leq 0$. But we also have
\[ \degree ( {\cal M} ) = -1 - \degree ( {\cal L} ) \geq 0. \]
Thus $\degree ( {\cal M} ) = 0$ and $\degree ( {\cal L} ) = -1$. By stability, these subsheaves must in fact correspond to line subbundles by Lemma \ref{injectivity}. Therefore, by Lemma \ref{lsubbsEth} there are at most four possibilities for each of the corresponding line subbundles of $E$ and $E^*$.
\par
Conversely, let $L \subset E$ and $M \subset E^*$ be line subbundles of degree $-1$ and $0$ respectively, and let $0 \to L \to F \to M \to 0$ be a nontrivial extension. We want to know when the composed map $g_{F} \colon F \to M \to E^*$ factorises via $W$. By Lemma \ref{extlift}, this is equivalent to $\delta(W)$ belonging to the kernel of the induced map
\[ g_{F}^{*} \colon H^{1}(X, \Sym^{2} E) \to H^{1}(X, \Homom(F, E)). \]
Clearly, this factorises via $H^{1}(X, \Homom(M, E)) = \C^3$.
\par
Now if $F$ were isomorphic to $E$, we would have a nonzero bundle map $E \to E^*$, contradicting genericity. Since $F$ and $E$ are stable, this implies that $h^{0}(X, \Homom(F,E)) = 0$. Therefore we have an exact cohomology sequence
\begin{multline*} 0 \to H^{0}(X, \Homom(L, E)) \to H^{1}(X, \Homom(M, E)) \to \\ H^{1}(X, \Homom(F, E)) \to H^{1}(X, \Homom(L, E)) \to 0. \end{multline*}
Since $h^{0}(X, \Homom(L, E)) = 1$, the map
\[ H^{1}(X, \Homom(M, E)) \to H^{1}(X, \Homom(F, E)) \]
has rank $2$. But we recall from the proof of Lemma \ref{genWstable} that the map $H^{1}(X, \Sym^{2} E) \to H^{1}(X, \Homom(M, E))$ is surjective. Therefore the rank of $g_{F}^*$ is also $2$, so its kernel is of dimension $4$.
\par
We now make a dimension count. There is a $0$-dimensional choice for $L$ and $M$, a $1$-dimensional choice for the extension $F$ and each $g_{F}^*$ has a $4$-dimensional kernel, giving a $\PP^3$ in $\PP^{5}_E$. Hence the set of extensions in $\PP^{5}_E$ admitting such an $F$ is of codimension at least $1$ in $\PP^{5}_E$. The proposition follows. \qed \\
\par
The last possibility is that $\dim(F|_{x} \cap E|_{x}) = 0$ for all but finitely many $x \in X$. Thus $F$ is an elementary transformation
\[ 0 \to {\cal F} \to {\cal E}^{*} \to T \to 0 \]
for some torsion sheaf $T$ of length $2$. The following result is in the spirit of Lange--Narasimhan \cite{LNar1983}, Prop.\ 1.1. We recall that a \textsl{2-secant} to a variety $Y \subset \PP^n$ is the projective linear span of two distinct points of $Y$ or else a line tangent to some point of $\psi ( \PP E)$.
\begin{lemma}
Let $W$ be a stable symplectic extension of $E^*$ by $E$. Then the number of degree $-1$ elementary transformations of $E^*$ which lift to isotropic subbundles of $W$ is bounded by the number of 2-secants to the surface $\psi(\PP E) \subset \PP^{5}_E$ which pass through $\langle \delta(W) \rangle$. \label{eltransecant}
\end{lemma}
\textbf{Proof}\\
Consider an elementary transformation $0 \to {\cal F} \to {\cal E}^{*} \to T \to 0$ such that $F$ lifts to an isotropic subbundle of $W$ of degree $-1$. Since $E$ and $E^*$ are stable and $\mu(E^{*}) > \mu(E)$, there are no maps $E^{*} \to E$. Therefore, by Thm.\ \ref{extprops} (ii) the sheaf ${\cal F}$ is the kernel of a symmetric $p \in \Prin(\Sym^{2}E)$ such that $\delta(W) = [p]$. There are three possibilities for $T$, which is identified with the image of $p \colon {\cal E}^{*} \to \sPrin(E)$:
\begin{enumerate}
\renewcommand{\labelenumi}{(\roman{enumi})}
\item $\C_{x} \oplus \C_{y}$ for distinct $x, y \in X$.
\item $\C_{x}^{\oplus 2}$ for some $x \in X$.
\item $\C_{2x}$ for some $x \in X$.
\end{enumerate}
We treat each case in turn.
\par
(i) Here $p$ is a sum of two principal parts $p_1$ and $p_2$ with simple poles along vectors $g \otimes g \in \Sym^{2}E|_x$ and $h \otimes h \in \Sym^{2}E|_y$ respectively. They are symmetric because $p$ is and decomposable because $p$ is not surjective to
\[ \frac{{\cal E}(x)}{{\cal E}} \quad \hbox{or} \quad \frac{{\cal E}(y)}{{\cal E}}. \]
We claim that $\psi \langle g \rangle \neq \psi \langle h \rangle$. For, otherwise there is a relation $[ p_{1} ] = \lambda [ p_{2} ]$ for some $\lambda \in \C^*$. But then $\delta(W)$ is proportional to $[ p_{1} ]$, say, so $\langle \delta(W) \rangle \in \psi(E)$. By Prop.\ \ref{strsemist}, then, $W$ is not stable, contrary to hypothesis.
\par
Thus $\langle \delta(W) \rangle$ lies on the secant spanned by the distinct points $\psi \langle g \rangle$ and $\psi \langle h \rangle$.\\
\par
(ii) Here $p$ is supported at one point $x \in X$. The image of the map $p \colon {\cal E}^{*} \to \sPrin(E)$ is isomorphic to $T$, so is of length $2$ and has only simple poles. Since $E$ is of rank $2$, the map $p$ is surjective to the subsheaf
\[ \frac{{\cal E}(x)}{{\cal E}} \subset \sPrin ( E ). \]
Therefore $p$ has a simple pole along some indecomposable vector in $\Sym^{2}E$ and is otherwise regular. Such a vector can be written as $g \otimes g + h \otimes h$ for some linearly independent $g, h \in E|_x$ (for example, because the image of $\langle g \rangle \mapsto \langle g \otimes g \rangle$ is nondegenerate in $\PP \: \Sym^{2}E|_x$). As above, if $\psi \langle g \rangle = \psi \langle h \rangle$ then $\langle \delta(W) \rangle$ belongs to the image of $\psi$, contradicting stability. Thus $\langle \delta(W) \rangle$ lies on the secant to $\psi \langle g \rangle$ and $\psi \langle h \rangle$.\\
\par
(iii) Again, $p$ is symmetric and supported at one point $x \in X$, where this time it has a double pole. Since $p$ is not surjective as a map ${\cal E}^{*} \to {\cal E}(2x)_x$, the double pole must be along a decomposable vector $g \otimes g \in \Sym^{2}E|_x$. The image of $p$ is a torsion sheaf of length $2$, so it must be equal to the subsheaf of ${\cal E}(2x) / {\cal E}$ of principal parts with poles of order up to $2$ in the direction of $g \subset E|_x$, and $p$ has poles in no direction other than $g \otimes g$.
\par
Recall that $\Upsilon \to \PP E$ is the line bundle $\pi^{*}K_{X} \otimes O_{\PP E}(2)$. Let $e, f$ be a local basis for $E$ near $x$ such that $e|_{x} = g$. Examining the local expression from the proof of Lemma \ref{psiemb} for a section of $\Upsilon$ vanishing to order 2 at $\langle g \rangle =:u $, we see that the contraction of such a section against a principal part of the form
\[ \frac{e \otimes e}{z^2} + \frac{\lambda e \otimes e}{z}, \]
for any $\lambda \in \C$, is regular. Thus $\langle [p] \rangle$ belongs to the embedded tangent space
\[ \PP \: \Ker \left( H^{0}( \PP E, \Upsilon )^{*} \to H^{0}( \PP E , \Upsilon - 2 \langle g \rangle )^{*} \right) \]
to $\psi ( \PP E )$ at $\psi \langle g \rangle$. In particular, $\langle [p] \rangle$ lies on a $2$-secant to $\psi ( \PP E )$. This completes the proof of Lemma \ref{eltransecant}.
 \qed \\
\\
\textbf{Remark:}  It is intriguing that in case (iii), in fact $\langle \delta(W) \rangle$ belongs to a particular line in the embedded tangent space to $\psi(\PP E)$, namely that spanned by classes of principal parts with single and double poles along $v \otimes v$. The exact relationship between secants to $\psi(\PP E)$ and subsheaves lifting to $W$ is still under investigation.\\
\\
We return to the proof of Theorem \ref{main}. By Lemma \ref{psiemb}, there are at most a finite number of points where $\psi$ fails to be an embedding. Therefore, the union of all 2-secants to $\psi ( \PP E )$ is of dimension at most 5 and through a general point of $\PP^{5}_E$ there pass a finite number of 2-secants to $\psi ( \PP E )$.
\par
Let $W$ be a general stable symplectic extension of $E^*$ by $E$. By Prop.\ \ref{FintEdrk1} and Lemma \ref{eltransecant}, the bundle $W$ has finitely many mutually nonisomorphic isotropic subbundles of rank $2$ and degree $-1$. Hence $W$ is represented in a finite number of fibres of $\mathbb{B}$.
\begin{prop}
Let $E \to X$ be a stable bundle of rank $2$ and degree $-1$ such that there are no nonzero maps $E \to E^*$. Let $W$ and $W^{\prime}$ be two extensions of $E^*$ by $E$. Then $W \cong W^{\prime}$ if and only if $\delta(W^{\prime}) = \lambda \delta(W)$ for some $\lambda \in \C^*$. \label{resphiinj}
\end{prop}
\textbf{Proof}\\
This is a special case of Narasimhan--Ramanan \cite{NR1969}, Lemma 3.3.\\
\par 
Write $U$ for the dense subset $\{ E \in \utmo : h^{0}(X, \Homom(E, E^{*})) = 0 \}$ of $\utmo$ and denote its complement $\Delta$. We observe that a general bundle $W \in \Mt$ cannot contain an isotropic subbundle $E$ belonging to $\Delta$: otherwise, the restriction of $\Phi$ to $\mathbb{B}|_{\Delta}$ would be dominant, which is impossible since $\dim ( \mathbb{B}|_{\Delta} ) = 9$ and $\dim ( \Mt) = 10$. Thus we can assume that the finitely many fibres of $\mathbb{B}$ in which our general $W$ is represented all lie over $U$. By Prop.\ \ref{resphiinj}, it occurs only once in each of these fibres. Thus the fibre of $\Phi$ over $W$ is finite. By semicontinuity, $\Phi$ is generically finite. Since $\mathbb{B}$ and $\Mt$ are of the same dimension, $\Phi$ is dominant. Since it a morphism of projective varieties, its image is closed, so it is actually surjective.

\subsection*{The degree of $\Phi$}
We use some results of Lange and Newstead to calculate the degree of $\Phi$.
\begin{thm}
Let $W$ be a general vector bundle of rank $n$ and degree $d$ over $X$. Suppose that $n \geq 4$ is even and $2d+4 \cong 0 \mod n$ with $\frac{2d+4}{n}$ odd. Then the number of subbundles of rank $2$ and maximal degree, counted with multiplicity, is equal to $\frac{n^3}{48} \left( n^{2} + 2 \right)$.
\end{thm}
\textbf{Proof}\\
Lange--Newstead \cite{LNew2002}, pp.\ 7--10.\\
\\
We check that a general $W \in \Mt$ is general enough in the sense of Lange and Newstead. We follow the notation of \cite{LNew2002}: for a subbundle $F \subseteq W$, we define
\[ s(W,F) = -\rank(W)\degree(F) + \rank(F)\degree(W) \]
and
\[ s_{n^{\prime}}(W) := \min \{ s(W,F) : F \subseteq W \hbox{ of rank $n^{\prime}$} \}. \]
For us, $n^{\prime}=2$, $g=2$ and $n=4$. We check the generality conditions, which are listed on p.\ 6 of \cite{LNew2002}:
\begin{enumerate}
\renewcommand{\labelenumi}{(\roman{enumi})}
\item $s_{n^{\prime}}(W) =  n^{\prime}(n - n^{\prime})(g-1)$.
\item $s_{n_1}(W) \geq  n_{1}(n - n_{1})(g-1)$ for all $n_{1} \in \{ 1, \ldots n^{\prime}-1 \}$.
\item $W$ has only finitely many maximal subbundles of rank $n^{\prime}$.
\end{enumerate}
\quad \\
For (i), note that $s_{2}(W) = \min \{ -4 \: \degree(F) : F \subseteq W \hbox{ of rank $2$} \} = 4$ since $W$ is stable but contains a Lagrangian subbundle of degree $-1$. On the other hand, $n^{\prime}(n - n^{\prime})(g-1) = 2(4-2)(2-1) = 4$.
\par
In (ii), the only value of $n_1$ that we need to check is $1$. We have
\begin{align*} s_{1}(W) &= \min \{ s(W,L) : L \hbox{ a line subb.\ of } W \} \\
&= \min \{ -4 \: \degree(L) : L \hbox{ a line subb.\ of } W \} \\
&= 4 \end{align*}
since $W$ is stable by Lemma \ref{genWstable}. On the other hand, we have
\[ n_{1}(n-n_{1})(g-1) = 3. \]
\newpage
As for (iii): we saw in the preceding section that apart from $E$, all maximal rank $2$ subbundles of a general $W$ lift from elementary transformations of $E^*$. By Lemma \ref{eltransecant}, there are at most finitely many isotropic subbundles of this type. It will suffice, therefore, to show that \textsl{all} subbundles of rank 2 and degree $-1$ of a general $W$ are isotropic.
\par
Let $F \subset W$ be such a subbundle. The symplectic form on $W$ restricts to a global section of $\bigwedge^{2}F^*$, that is, $(\det F)^{-1}$. To compute this line bundle, we note that ${\cal F}$ fits into a short exact sequence
\[ 0 \to {\cal F} \to {\cal E}^{*} \xrightarrow{p} T \to 0 \]
where $T$ is a torsion sheaf of length $2$ and $p \in \Prin(\Homom(E^{*}, E))$ is such that $\delta(W) = [p]$. We have $\Supp(T) = \Supp(p) = \{ x, y \}$ for points $x, y \in X$ which are not necessarily distinct. Thus
\[ (\det F)^{-1} = (\det E)(x+y). \]
This has a nonzero section only if it is effective. In order for $W$ to contain a nonisotropic subbundle of rank 2 and degree $-1$, then, $\langle \delta(W) \rangle$ must lie on a line joining points in the fibres of $\psi ( \PP E )$ over $x$ and $y$ for some $x+y$ such that $\det (E)(x+y)$ is effective. For each such divisor $x+y$ the family of such classes $\langle \delta(W) \rangle$ in $\PP^{5}_E$ has dimension at most
\[ \dim ( \PP E|_{x} ) + \dim ( \PP E|_{y} ) + \dim ( \PP^{1} ) = 3. \]
By genericity, $(\det E)^{-1}$ is not effective so $\det E$ is of the form $\strs(q - r - s)$ with $\strs(r+s) \neq K_X$. We determine the divisors $x + y$ such that the bundle $(\det E)(x+y)$ is effective, that is,
\[ O_{X}(q-r-s+x+y) = O_{X}(t) \hbox{ for some } t \in X. \]
Equivalently, we seek all solutions $x, y, t \in X$ to the equation
\[ O_{X}(q+x+y) = O_{X}(r+s+t). \]
\par
Let $t \in X$. If $|r+s+t|$ has a base point then this must be $r$ or $s$ because $r+s \not\in |K_{X}|$, so in particular $q$ is not a base point. Therefore there is a unique divisor $x+y$ such that $q+x+y \sim r+s+t$. Letting $t$ vary, we find a one-dimensional family of pairs $x+y$ such that $\det (E)(x+y)$ is effective.
\par
Furthermore, if we had started with a different representative $q^{\prime} - r^{\prime} - s^{\prime}$ for $\det E$ then we would have obtained the same family of $x+y$. Explicitly,
\begin{align*} O_{X}(q+x+y) = O_{X}(r+s+t) & \Leftrightarrow q-r-s \sim t-x-y \\
& \Leftrightarrow q^{\prime} - r^{\prime} - s^{\prime} \sim t-x-y \hbox{ by hypothesis}\\
& \Leftrightarrow O_{X}(q^{\prime}+x+y) = O_{X}(r^{\prime}+s^{\prime}+t). \end{align*}
Also, since for each $t \in X$ there exists unique such $x+y$, these are all the solutions.
\par
Putting all this together, we find that the locus of $\langle \delta(W) \rangle$ in $\PP^{5}_E$ such that $W$ can contain a nonisotropic subbundle of rank 2 and degree $-1$ is of dimension at most 4. Therefore, for general $\langle \delta(W) \rangle \in \PP^{5}_E$, all maximal subbundles lifting from elementary transformations of ${\cal E}^*$ are isotropic.\\
\\
Thus $W$ has finitely many maximal subbundles of rank $2$. In summary, we find that a general symplectic bundle is indeed general in the sense of Lange and Newstead.\\
\\
In our case, the number $\frac{n^3}{48} \left( n^{2} + 2 \right)$ is $24$. To verify that $\deg \Phi = 24$, we need to check that all the subbundles of a general $W \in \Mt$ are isotropic and distinct. We have just checked isotropy. For distinctness: let $E \subset W$ be a maximal subbundle of rank $2$ and suppose that $h^{0}(X, \Homom(E, W)) \geq 2$. By generality, some copy of $E$ in $W$ is isotropic, so $W/E \cong E^*$ and we get a diagram
\[ \xymatrix{ 0 \ar[r] & E \ar[r] & W \ar[r] & E^{*} \ar[r] & 0 \\ & & E \ar[u] \ar[ur] & & } \]
where the composed map is nonzero. But this yields a nonzero map $E \to E^*$, contradicting genericity.\\
\\
In summary, $\Phi \colon \mathbb{B} \to \Mt$ is generically finite of degree $24$. This completes the proof of Theorem \ref{main}, our main result. \qed\\
\\
\textbf{Remark:} Some fibres of $\Phi$ may be of positive dimension. (This is analogous, for example, to the fact that there exist vector bundles of rank $2$ and degree $0$ over a curve of genus $3$ which admit infinitely many maximal line subbundles.) If, for example, $\psi \colon \PP E \to \PP^{5}_E$ failed to be an embedding over a locus of positive dimension in $\PP E$ then there might be points of $\PP^{5}_E$ lying on infinitely many $2$-secants to $\psi(\PP E)$. Then Lemma \ref{eltransecant} might lead us to expect that that infinitely many degree $-1$ elementary transformations of $E$ would lift to the corresponding extensions of $E^*$ by $E$. This will be a subject of future study.

\section{Applications and future work}

We give one immediate application of this construction.
\begin{cor}
The moduli space $\Mt$ is of Kodaira dimension $-\infty$.
\label{KDMt}
\end{cor}
\textbf{Proof}\\
Recall that we have a proper surjective map $\mathbb{B} \to {\cal U}_X(2,-1)$ whose fibres are projective spaces. Hence $\mathbb{B}$ is uniruled, so is of Kodaira dimension $-\infty$. Since by Thm.\ \ref{main} there is a morphism $\mathbb{B} \to \Mt$ which is surjective and generically finite, the same is true for $\Mt$.  \qed \\

Another application of Thm.\ \ref{main} is given, in the last chapter of \cite{Hit2005a}, to the study of theta divisors of symplectic vector bundles over a curve of genus $2$. We sketch the result. Let $W$ be a semistable symplectic bundle of rank $4$ and consider the set
\[ S(W) = \left\{ L \in J^{1}_{X} : h^{0}(X, \Homom(L, W)) > 0 \right\}. \]
For general $W$, this set is the support of an even $4\Theta$ divisor $D(W)$, and this defines a rational map $D \colon \Mt \dashrightarrow |4\Theta|_{+} = \PP^9$. The line bundle $D^{*}O(1) =: \Xi$ is called the \textsl{determinant} bundle on $\Mt$.
\par
For certain $W$, however, $S(W)$ is the whole Jacobian; equivalently, $W$ is a base point of $|\Xi|$. If one wishes to find such bundles, by Thm.\ \ref{main} it suffices to look for them in the spaces $\PP^{5}_E$. In \cite{Hit2005a}, Chap.\ 6, we find necessary conditions on $E \in \utmo$ for the existence of an extension $W$ in $\PP^{5}_E$ such that $S(W) = J^{1}_X$. We hope that this will lead to a full description of the base locus of $|\Xi|$. See Beauville \cite{Bea2005}, \cite{Bea1995} and Raynaud \cite{Ray1982} for general results on this type of question.\\
\par
It is natural to ask whether a similar description of $\Mt$ exists in higher genus. One calculates that for the dimensions to work out as in the genus $2$ case, we should consider symplectic extensions $0 \to E \to W \to E^{*} \to 0$ where $E \to X$ is a vector bundle of rank $2$ and degree $1-g$. The $\PP^{5}_E$ are replaced by projective spaces of dimension $6g-7$. However, there is no Poincar\'e bundle over ${\cal U}_{X}(2, g-1) \times X$ if $g-1$ is not relatively prime to $2$, so we are forced to restrict to the case of even genus.
\par
If symplectic bundles in higher even genus are still general enough then, by Lange--Newstead \cite{LNew2002}, Prop.\ 2.4, the problem involves determining those elementary transformations $0 \to {\cal F} \to {\cal E}^{*} \to T \to 0$, where $T$ is a torsion sheaf of length $2g-2$, which lift to the above $W$. Now one expects that the $(2g-2)$nd secant variety to a nondegenerate surface in $\PP^{6g-7}$ with a finite number of singular points will be of dimension $6g-7$, so Lemma \ref{eltransecant} might lead us to conjecture that such an extension would indeed have only finitely many Lagrangian subbundles of degree $g-1$. This will be studied in the future.

\noindent Institut f\"ur Algebraische Geometrie,\\
Universit\"at Hannover,\\
Welfengarten 1,\\
30167 Hannover,\\
Germany.\\
Tel.: +49 (0)511 762 3206\\
Fax.: +49 (0)511 762 5803\\
\texttt{hitching@math.uni-hannover.de}


\begin{thebibliography}{99}

\bibitem {Ati1957} Atiyah, M.\ F.: \textsl{Complex analytic connections in fibre bundles}, Trans.\ Amer.\ Math.\ Soc.\ \textbf{85} (1957), pp.\ 181--207.

\bibitem {Bea2005} Beauville, A.: \textsl{Vector bundles on curves and theta functions}, arXiv AG/0502179.

\bibitem {Bea2004} Beauville, A.: \textsl{Vector bundles and theta functions on curves of genus 2 and 3}, arXiv AG/0406030, to appear in Amer.\ J.\ of Math.

\bibitem {Bea1995} Beauville, A.: \textsl{Vector Bundles on Curves and Generalized Theta Functions: Recent Results and Open Problems}, MSRI Series, Vol.\ 28, 1995, pp.\ 17--33.

\bibitem {BL2004} Birkenhake, Ch.; Lange, H.: \textsl{Complex Abelian varieties}, second edition. Grundlehren der Mathematischen Wissenschaften 302. Springer-Verlag, Berlin, 2004.

\bibitem {GH1978} Griffiths, P.; Harris J.: \textsl{Principles of Algebraic Geometry}, USA, John Wiley and Sons Inc., 1978.

\bibitem {Har1977} Hartshorne, R.: \textsl{Algebraic Geometry}, Graduate Texts in Mathematics 52, New York, Springer-Verlag Inc., 1977.

\bibitem {Hit2005a} Hitching, G.\ H.: \textsl{Moduli of symplectic bundles over curves}, doctoral thesis, University of Durham, 2005.

\bibitem {Hit2005b} Hitching, G.\ H.: \textsl{A remark on subbundles of symplectic and orthogonal vector bundles over curves}, arXiv AG/0407323 (v3), 2005, revised version submitted.

\bibitem {Kem1983} Kempf, G.: \textsl{Abelian Integrals}, Monograf\'{\i}as del Instituto de Matem\'aticas 13, Universidad Nacional Aut\'onoma de M\'exico, Mexico, 1983.

\bibitem {KS1988} Kempf, G.; Schreyer, F.-O.: \textsl{A Torelli theorem for osculating cones to the theta divisor}, Compositio Math.\ \textbf{67} (1988), no.\ 3, p.\ 343--353.

\bibitem {Lan1993} Lang, S.: \textsl{Algebra}, USA, Addison--Wesley Publishing Company, 1993.

\bibitem {LNar1983} Lange, H.; Narasimhan, M.\ S.: \textsl{Maximal subbundles of rank two vector bundles on curves}, Math.\ Ann.\ \textbf{266} (1983), no.\ 1, pp.\ 55--72.

\bibitem {LNew2002} Lange, H.; Newstead, P.\ E.: \textsl{Maximal subbundles and Gromov--Witten invariants}. A tribute to C.\ S.\ Seshadri (Chennai, 2002), pp.\ 310--322, Trends Math., Birkh\"auser, Basel, 2003.

\bibitem  {LeP1997} Le Potier, J.: \textsl{Lectures on vector bundles} (translated by A.\ Maciocia), Cambridge Studies in Advanced Mathematics, 54, Cambridge University Press, Cambridge, 1997.

\bibitem {Mat2002} Matsuki, K.: \textsl{Introduction to the Mori program}, Universitext, Springer-Verlag, New York, 2002.

\bibitem {NR1987} Narasimhan, M.\ S.; Ramanan, S.: \textsl{$2\Theta$ linear systems on Abelian varieties}, \textsl{Vector bundles on algebraic varieties} (Bombay, 1984), pp.\ 415--427, Tata Inst.\ Fund.\ Res.\ Stud.\ Math., 11, Tata Inst.\ Fund.\ Res., Bombay, 1987.

\bibitem {NR1969} Narasimhan, M.\ S.; Ramanan, S.: \textsl{Moduli of vector bundles on a compact Riemann surface}, Ann.\ Math.\ (2) \textbf{89} (1969), pp. 14--51.

\bibitem {Pau2002} Pauly, C.: \textsl{Self-duality of Coble's quartic hypersurface and applications}, Michigan Math.\ J.\ \textbf{50} (2002), no.\ 3, pp.\ 551--574.

\bibitem {Ram1981} Ramanan, S.: \textsl{Orthogonal and spin bundles over hyperelliptic curves}, \textsl{Geometry and Analysis: Papers dedicated to the memory of V. K. Patodi}, Berlin, Springer--Verlag, 1981.

\bibitem {Rthn1996a} Ramanathan, A.: \textsl{Moduli for principal bundles over algebraic curves I}, Proc.\ Indian Acad.\ Sci.\ Math.\ Sci.\ \textbf{106} (1996), no.\ 3, pp.\ 301--328.

\bibitem {Rthn1996b} Ramanathan, A.: \textsl{Moduli for principal bundles over algebraic curves II}, Proc.\ Indian Acad.\ Sci.\ Math.\ Sci.\ \textbf{106} (1996),  no.\ 4, pp.\ 421--449.

\bibitem {Rthn1975} Ramanathan, A.: \textsl{Stable principal bundles on a compact Riemann surface}, Math.\ Ann.\ 213 (1975), pp.\ 129--152.

\bibitem {Ray1982} Raynaud, M.: \textsl{Sections des fibr\'es vectoriels sur une courbe}, Bull.\ Soc.\ Math.\ France \textbf{110} (1982), no.\ 1, pp.\ 103--125.

\bibitem {Ses1982} Seshadri, C.\ S.: \textsl{Fibr\}es Vectoriels sur les courbes alg\'ebriques} (notes by J.\ M.\ Drezet), Ast\'erisque 96, Soc.\ Math.\ de France, Paris, 1982.

\end{thebibliography}
\end{document}